\documentclass[reqno,11pt]{amsart}
\usepackage{amsmath, latexsym, amsfonts, amssymb, amsthm, amscd}

\numberwithin{equation}{section}

\newtheorem{theorem}{Theorem}[section]
\newtheorem{lemma}[theorem]{Lemma}
\newtheorem{prop}[theorem]{Proposition}

\newtheorem{definition}[theorem]{Definition}

\newcommand{\ga}{\alpha}

\newcommand{\gga}{\gamma}            

\newcommand{\gep}{\varepsilon}       

\newcommand{\gl}{\lambda}

\newcommand{\go}{\omega}

\newcommand{\cF}{{\ensuremath{\mathcal F}} }

\newcommand{\cE}{{\ensuremath{\mathcal E}} }

\newcommand{\cN}{{\ensuremath{\mathcal N}} }

\newcommand{\cL}{{\ensuremath{\mathcal L}} }

\newcommand{\cZ}{{\ensuremath{\mathcal Z}} }

\newcommand{\bbC}{{\ensuremath{\mathbb C}} }

\newcommand{\bbE}{{\ensuremath{\mathbb E}} }
\newcommand{\E}{{\ensuremath{\mathbb E}} }

\newcommand{\bbN}{{\ensuremath{\mathbb N}} }
\newcommand{\N}{{\ensuremath{\mathbb N}} }

\newcommand{\bbP}{{\ensuremath{\mathbb P}} }

\newcommand{\bbR}{{\ensuremath{\mathbb R}} }
\newcommand{\R}{{\ensuremath{\mathbb R}} }

\title{A conservative evolution of the Brownian excursion}

\author{Lorenzo Zambotti}
\address{Laboratoire de Probabilit{\'e}s et Mod\`eles Al\'eatoires (CNRS U.M.R. 7599) \\ Universit{\'e} Paris 6
-- Pierre et Marie Curie, U.F.R. Mathematiques, Case 188, 4 place
Jussieu, 75252 Paris cedex 05, France }
\email{zambotti\@@ccr.jussieu.fr}

\date{}

\begin{document}

\begin{abstract}
We consider the problem of conditioning the Brownian excursion to
have a fixed time average over the interval $[0,1]$ and we study an
associated stochastic partial differential equation with reflection
at 0 and with the constraint of conservation of the space average.
The equation is driven by the derivative in space of a space-time
white noise and contains a double Laplacian in the drift. Due to the
lack of the maximum principle for the double Laplacian, the standard
techniques based on the penalization method do not yield existence
of a solution.
\\ \\
2000 \textit{Mathematics Subject Classification: 60J65, 60G15,
60H15, 60H07, 37L40}
\\ \\
\noindent\textit{Keywords: Brownian meander; Brownian excursion;
singular conditioning; Stochastic partial differential equations
with reflection}
\end{abstract}

\maketitle

\section{Introduction}

The aim of this paper is to construct a stochastic evolution, whose
invariant measure is the law of the Brownian excursion $(e_\theta,
\theta\in[0,1])$ conditioned to have a fixed average $\int_0^1
e_\theta \, d\theta = c>0$ over the interval $[0,1]$.

Since the distribution of the random variable $\int_0^1 e_\theta \,
d\theta$ is non-atomic, and the Brownian excursion is not a Gaussian
process, it is already not obvious that such conditioning is well
defined. The first part of the paper will be dedicated to this
problem: we shall write down the density of the random variable
$\int_0^1 e_\theta \, d\theta$ and a regular conditional
distribution of the law of $(e_\theta, \theta\in[0,1])$, given
$\int_0^1 e_\theta \, d\theta$. The same will be done for the
Brownian meander $(m_\theta, \theta\in[0,1])$.

After this is done, we shall turn to the problem of constructing a
natural stochastic dynamics associated with the constructed
conditioned laws. We recall that a stochastic dynamics whose
invariant measure is the law of the Brownian excursion has been
studied in \cite{nupa} and \cite{za}, where a stochastic partial
differential equation with reflection and driven by space-time white
noise is proven to be well posed and associated with a Dirichlet
form with reference measure given by the law of the Brownian
excursion.

In the present case, we shall see that a natural dynamics with the
desired properties solves a stochastic partial differential equation
of the fourth order with reflection and driven by the derivative in
space of a space-time white noise:
\begin{equation}\label{00}
\left\{ \begin{array}{ll} {\displaystyle \frac{\partial u}{\partial
t}=- \frac{\partial^2 }{\partial \theta^2} \left(\frac{\partial^2
u}{\partial \theta^2} + \eta \right)  + \sqrt 2 \,
\frac{\partial}{\partial\theta} \dot{W}, }
\\ \\
{\displaystyle u(t,0)=u(t,1)= \frac{\partial^3 u}{\partial
\theta^3}(t,0)=\frac{\partial^3 u}{\partial \theta^3}(t,1)=0 }
\\ \\
u(0,\theta)=x(\theta)
\end{array} \right.
\end{equation}
where $\dot W$ is a space-time white noise on
$[0,+\infty)\times[0,1]$, $u$ is a continuous function of
$(t,\theta)\in [0,+\infty)\times[0,1]$, $\eta$ is a locally finite
positive measure on $(0,+\infty)\times[0,1]$, subject to the
constraint:
\begin{equation}\label{contact}
u\geq 0, \qquad \int_{(0,+\infty)\times[0,1]} u \, d\eta \, = \, 0.
\end{equation}
This kind of equations arise as scaling limit of fluctuations of
conservative interface models on a wall, as shown in \cite{f}, where
however different boundary conditions are considered.

Indeed, notice that the boundary conditions in \eqref{00} are mixed,
i.e. Dirichlet for $u$ and Neumann for
$\frac{\partial^2u}{\partial\theta^2}$. In \cite{deza} and \cite{f}
a similar equation, with Neumann boundary conditions for $u$ and
$\frac{\partial^2u}{\partial\theta^2}$, has been studied, together
with the scaling limit of interface models mentioned above. In that
case it was possible to prove pathwise uniqueness and existence of
strong solutions for the SPDE. In the case of \eqref{00} we can only
prove existence of weak solutions, since we have failed to obtain a
uniqueness result: see subsection \ref{remuni} below.

The dynamics is anyway uniquely determined by a natural
infinite-dimensional Dirichlet form on the path space of the
Brownian excursion, to which it is associated. See Theorem
\ref{main?} below.

We also notice that the Brownian meander $(m_\theta,
\theta\in[0,1])$ conditioned to have a fixed average and the density
of $\int_0^1 m_\theta \, d\theta$ appear in an infinite-dimensional
integration by parts formula in \cite[Corollary 6.2]{deza}.

\section{The main results}\label{tmr}

In this section we want to present the setting and the main results
of this paper.

\subsection{Conditioning the Brownian excursion to have a fixed time average}
\label{exx}

Let $(e_t,t\in[0,1])$ be the normalized Brownian excursion, see
\cite{reyo}, and $(\beta_t,t\in[0,1])$ a Brownian bridge between $0$
and $0$. Let $\{m,\hat m,b\}$ be a triple of processes such that:
\begin{enumerate}
\item $m$ and $\hat m$ are independent copies of a Brownian meander on
$[0,1]$
\item conditionally on $\{m,\hat m\}$, $b$ is a Brownian bridge on
$[1/3,2/3]$ from $\frac 1{\sqrt 3}\, m_1$ to $\frac 1{\sqrt 3}\,
\hat m_1$
\end{enumerate}
We introduce the continuous processes:
\begin{align}
\label{v1} & v_t \, := \, \left\{ \begin{array}{ll} {\displaystyle
\frac 1{\sqrt 3}\, m_{3t}, \qquad \quad t\in[0,1/3] }
\\ \\
{\displaystyle b_t, \qquad \qquad t\in[1/3,2/3],}
\\ \\
{\displaystyle \frac 1{\sqrt 3}\, \hat m_{1-3t}, \quad t\in[2/3,1],}
\end{array} \right.
\\
\label{v2} & V^c_t \, := \, \left\{ \begin{array}{ll} {\displaystyle
v_t, \qquad \quad t\in[0,1/3]\cup[2/3,1] }
\\ \\
{\displaystyle v_t + 18\left(9\, t\, (1-t)-2\right)\left(c-\int_0^1
v\right), \quad t\in[1/3,2/3].}
\end{array} \right.
\end{align}
Notice that $\int_0^1 V^c_t \, dt = c$. We set for all $\go\in
C([0,1])$:
\[
\rho^c(\go) \,:= \, \exp\left\{ -162\left(\int_0^\frac 13
\left(\go_r+\go_{1-r}\right) dr + \frac{\go_\frac 13 +\go_\frac
23}6- c\right)^2 - \frac 32 (\go_{\frac 23}-\go_{\frac 13})^2
\right\}.
\]
The first result of this paper is
\begin{theorem}\label{4.1} Setting for all $c\geq 0$
\[
p_{\langle e,1\rangle}(c) \, =: \, 27\sqrt{\frac6{\pi^3}} \,
\bbE\left[\rho^c(V^c) \, 1_{\{V^c_t\geq 0, \ \forall t\in[0,1]\}}
\right],
\]
and for all bounded Borel $\Phi:C([0,1])\mapsto\bbR$ and $c>0$
\[
\bbE\left[\Phi(e) \, | \, \langle e,1\rangle =c \right] \, := \,
\frac 1{\cZ_c} \ \bbE\left[\Phi\left(V^c\right)\, \rho^c(V^c) \, \,
1_{\{V^c_t\geq 0, \ \forall t\in[0,1]\}} \right],
\]
where $\cZ_c>0$ is a normalization factor, we have
\begin{enumerate}
\item $p_{\langle e,1\rangle}$ is the density of $\langle e,1\rangle$
on $[0,\infty)$, i.e.
\[
\bbP(\langle e,1\rangle\in dc)\, = \, p_{\langle e,1\rangle}(c)
\,1_{\{c\geq 0\}} \, dc.
\]
Moreover $p_{\langle e,1\rangle}$ is continuous on $[0,\infty)$,
$p_{\langle e,1\rangle}(c)>0$ for all $c\in(0,\infty)$ and
$p_{\langle e,1\rangle}(0)=0$.
\item $(\bbP\left[e\in\cdot \, | \, \langle e,1\rangle =c \right], c>0)$
is a regular conditional distribution of $e$ given $\langle
e,1\rangle$, i.e.
\[
\bbP(e\in\cdot\, , \langle e,1\rangle\in dc ) \, = \,
\bbP\left[e\in\cdot \, | \, \langle e,1\rangle =c \right] \
p_{\langle e,1\rangle}(c) \,1_{\{c> 0\}} \, dc.
\]
\end{enumerate}
\end{theorem}
\medskip\noindent
In section \ref{meander} below we state and prove analogous results
for the Brownian meander.

\subsection{Two Hilbert spaces}

For the study of the stochastic partial differential equation
\eqref{00} we need to introduce some notation. We denote by $\langle
\cdot,\cdot\rangle_L$ the canonical scalar product in $L^2(0,1)$:
\[
\langle h,k \rangle_L \, := \, \int_0^1 h_\theta \, k_\theta \,
d\theta.
\]
We denote by $A$ the realization in $L^{2}(0,1)$ of
$\partial^2_\theta$ with Neumann boundary condition at $0$ and $1$,
i.e.:
\begin{equation}\label{A}
D(A) \, := \, \{h\in H^2(0,1): \, h'(0)=h'(1)=0\}, \qquad A \, := \,
\frac{\partial^2}{\partial\theta^2}.
\end{equation}
Notice that $A$ is self-adjoint in $L^2(0,1)$. We introduce a
notation for the {\em average} of $h\in L^2(0,1)$:
\[
\overline h \, := \, \int_0^1 h \, = \, \langle h,1\rangle_L.
\]
Then we also set for all $c\in\bbR$:
\[
L^2_c \, := \, \left\{h\in L^2(0,1): \ \overline h = c\right\}
\]
Now we define the operator $Q:L^2(0,1)\mapsto L^2(0,1)$:
\begin{align*}
Qh(\theta) & \, := \, \int_0^1 q(\theta,\sigma) \, h_\sigma \,
d\sigma, \qquad {\rm where:}
\\
q(\theta,\sigma) & \, := \, \theta\wedge \sigma +
\frac{\theta^2+\sigma^2}2 - \theta-\sigma+\frac 43, \qquad
\theta,\sigma\in[0,1].
\end{align*}
Then a direct computation shows that for all $h\in L^2(0,1)$:
\[
\langle Qh,1\rangle_L \, = \, \langle h,1\rangle_L, \qquad -AQh \, =
\, h-\overline h,
\]
i.e. $Q$ is the inverse of $-A$ on $L^2_0$ and conserves the
average. Then we define $H$ as the completion of $L^2(0,1)$ with
respect to the scalar product:
\[
(h,k)_H \, := \, \langle Q h,k\rangle_L.
\]
For all $c\in\bbR$ we also set:
\[
H_c \, := \, \left\{ h\in H: \ (h,1)_H = c \right\}.
\]
We remark that $H$ is naturally interpreted as a space of
distributions, in particular as the dual space of $H^1(0,1)$.
%
%

\smallskip\noindent
We also need a notation for the realization $A_D$ in $L^{2}(0,1)$ of
$\partial^2_\theta$ with Dirichlet boundary condition at $0$ and
$1$, i.e.:
\begin{equation}\label{A_D}
D(A_D) \, := \, \{h\in H^2(0,1): \, h(0)=h(1)=0\}, \qquad A_D \, :=
\, \frac{\partial^2}{\partial\theta^2}.
\end{equation}
Notice that $A_D$ is self-adjoint and invertible in $L^2(0,1)$, with
inverse:
\begin{equation}\label{Q_D}
Q_Dh(\theta) \, = \, (-A_D)^{-1}h(\theta) \, := \, \int_0^1
\left(\theta\wedge \sigma-\theta\sigma\right) \, h_\sigma \,
d\sigma, \qquad \theta\in[0,1].
\end{equation}

\subsection{Weak solutions of (\ref{00})}
We state now the precise meaning of a solution to (\ref{00}).
\begin{definition}
\label{solution} Let $u_0\in C([0,1])$, $u_0\geq 0$, $\int_0^1
u_0>0$, $u_0(0)=u_0(1)=0$. We say that $(u,\eta,W)$, defined on a
filtered complete probability space $(\Omega,\bbP,\cF,\cF_t)$, is a
{\em weak} solution to (\ref{00}) on $[0,T]$ if
\begin{enumerate}
\item a.s. $u\in C((0,T]\times [0,1])$, $u\geq 0$ and $u\in C([0,T];H)$
\item a.s. $u_t(0)=u_t(1)=0$ for all $t\geq 0$
\item a.s. $\eta$ is a positive measure on $(0,T]\times (0,1)$, such that $\eta([\delta,T]\times[\delta,1-\delta])<\infty$
for all $\delta>0$
\item $(W(t,\theta))$ is a Brownian sheet, i.e. a centered Gaussian process such that
\[
\bbE\, [W(t,\theta) \, W(t',\theta')]=t\wedge t' \cdot
\theta\wedge\theta', \qquad t,t'\geq 0, \ \theta,\theta'\in[0,1]
\]
\item $u_0$ and $W$ are independent and the process $t\mapsto (u_t(\theta), 
W(t,\theta))$ is $(\cF_t)$-adapted for all $\theta\in[0,1]$ and $I$
interval in $[0,1]$
\item for all $h\in C^4([0,1])$ such that $h'(0)=h'(1)=h''(0)=h''(1)=0$ and for all $0<\delta\leq t\leq T$:
\begin{eqnarray}\nonumber
\langle u_t,h \rangle_L & = & \langle u_\delta,h \rangle_L -
\int_\delta^t \langle u_s, A_D Ah \rangle_L \, ds
\\ \nonumber \\ \label{distr}& & - \int_\delta^t \int_0^1 
Ah_\theta \, \eta(ds,d\theta) - \sqrt 2  \int_\delta^t \int_0^1
h'_\theta \, W(ds,d\theta)
\end{eqnarray}
\item a.s. the {\em contact property} holds: {\rm supp}$(\eta)\subset\{(t,\theta):u_t(\theta)=0\}$, i.e.
\[
\int_{(0,T]\times [0,1]} u \, d\eta=0.
\]
\end{enumerate}
\end{definition}
\noindent

\subsection{Function spaces}
Notice that for all $c\in\bbR$, $H_c=c1+H_0$ is a closed affine
subspace of $H$ isomorphic to the Hilbert space $H_0$. If $J$ is a
closed affine subspace of $H$ space, we denote by $C_b(J)$,
respectively $C_b^1(J)$, the space of all bounded continuous
functions on $J$, resp. bounded and continuous together with the
first Fr\'echet derivative. We also denote by ${\rm Lip}(J)$ the set
of all $\varphi\in C_b(J)$ such that:
\[
[\varphi]_{{\rm Lip}(J)} \, := \, \sup_{h\ne k} \,
\frac{|\varphi(h)-\varphi(k)|}{\|h-k\|_J} \, < \, \infty.
\]
Finally, we define ${\rm Exp}_A(H)\subset C_b(H)$ as the linear span
of $\{\cos((h, \cdot)_H), \sin((h,\cdot)_H): h\in D(A_D\, A)\}$.

To $\varphi\in C_b^1(H_c)$ we associate a gradient
$\nabla_{H_0}\varphi:H_c\mapsto H_0$, defined by:
\begin{equation}\label{gradiente}
\left.\frac d{d\varepsilon} \, \varphi(k+\gep \, h) \right|_{\gep=0}
\, = \, (\nabla_{H_0}\varphi(k),h)_{H}, \qquad \forall \ k\in A, \
h\in H_0.
\end{equation}
The important point here is that we only allow derivatives along
vectors in $H_0$ and the gradient is correspondingly in $H_0$. In
particular, by the definition of the scalar product in $H$, each
$\varphi\in {\rm Exp}_A(H)$ is also Fr\'echet differentiable in the
norm of $L^2(0,1)$; then, denoting by $\nabla \varphi$ the gradient
in the Hilbert structure of $L^2(0,1)$, we have
\begin{equation}\label{gradients}
\nabla_{H_0}\varphi \, = \, (-A)\nabla\varphi, \qquad \forall \
\varphi\in {\rm Exp}_A(H).
\end{equation}

\subsection{The stochastic dynamics}
We are going to state the result concerning equation \eqref{00}. We
denote by $X_t:H^{[0,\infty[}\mapsto H$ the coordinate process and
we define
\[
\nu_c := \bbP\left[e\in\cdot \, | \, \langle e,1\rangle =c \right] =
\text{law of } e \text{ conditioned to have average  }c,
\]
which is well defined by Theorem \ref{4.1}. We notice that the
support of $\nu_c$ in $H$ is
\[
K_c := \text{closure in $H$ of } \left\{h\in L^2(0,1): h\geq 0, \
\langle h,1\rangle_L =1 \right\},
\]
and the closed affine hull in $H$ of $K_c$ is $H_c$.

\medskip
Then the second result of this paper is
\begin{theorem}\label{main?} Let $c>0$.
\begin{itemize}
\item[(a)]
The bilinear form ${\cE}={\cE}_{\nu_c,\|\cdot\|_{H_0}}$ given by
\[
{\cE}(u,v) \, := \, \int_{K_c} (\nabla_{H_0} u,\nabla_{H_0} v)_{H}
\, d\nu_c, \qquad u,\,v\in C^1_b(H_c),
\]
is closable in $L^2(\nu_c)$ and its closure $(\cE,D(\cE))$ is a
symmetric Dirichlet Form. Furthermore, the associated semigroup
$(P_t)_{t\geq 0}$ in $L^2(\nu_c)$ maps $L^\infty(\nu_c)$ in
$C_b(K_c)$.
\item[(b)] For any $u_0=x\in K_c\cap C([0,1])$ there exists a weak
solution $(u,\eta,W)$ of \eqref{00} such that the law of $u$ is
$\bbP_x$.
\item[(c)] $\nu_c$ is invariant for $(P_t)$, i.e.
$\nu_c(P_tf)=\nu_c(f)$ for all $f\in C_b(K_c)$ and $t\geq 0$.
\end{itemize}
\end{theorem}
\noindent By Theorem \ref{main?}, we have a Markov process which
solves \eqref{00} weakly and whose invariant measure is the law of
$e$ conditioned to have average equal to $c$.

\subsection{Remarks on uniqueness of solutions to \eqref{00}} \label{remuni}
We expect equation \eqref{00} to have pathwise-unique solutions,
since this is typically the case for monotone gradient systems: this
is always true in finite dimensions, see \cite{cepa}, and has been
proven in several interesting infinite-dimensional situations, see
\cite{nupa} and \cite{deza}. In the present situation, the
difficulty we encountered in the proof of uniqueness of \eqref{00}
is the following: because of the boundary condition
$u(t,0)=u(t,1)=0$ and of the reflection at $0$, it is expected that
the reflecting measure $\eta$ has infinite mass on
$]0,T]\times[0,1]$; this is indeed true for second order SPDEs with
reflection: see \cite{za4}. If this is the case, then it becomes
necessary to localize in $]0,1[$ in order to prove a priori
estimates; however, in doing so one loses the crucial property that
the average is constant. In short, we were not able to overcome
these two problems.

\section{Conditioning $e$ on its average}

\subsection{An absolute continuity formula} \label{secabco}

Let $(X_t)_{t\in[0,1]}$ be a continuous centered Gaussian process
with covariance function $q_{t,s}:=\bbE[X_t\,X_s]$. We have in mind
the case of $X$ being a Brownian motion or a Brownian bridge. In
this section we consider two processes $Y$ and $Z$, both defined by
linear transformations of $X$, and we write an absolute continuity
formula between the laws of $Y$ and $Z$.

For all $h$ in the space $M([0,1])$ of all signed measures with
finite total variation on $[0,1]$ we set:
\[
Q:M([0,1])\mapsto C([0,1]), \qquad Q\gl(t) \, := \, \int_0^1 q_{t,s}
\, \gl(ds), \quad t\in[0,1].
\]
We denote by $\langle \cdot,\cdot \rangle:C([0,1])\times
M([0,1])\mapsto\bbR$ the canonical pairing,
\[
\langle h,\mu \rangle \, := \, \int_0^1 h_t \, \mu(dt).
\]
where a continuous function $k\in C([0,1])$ is identified with $k_t
\,dt\in M([0,1])$. We consider $\gl,\mu\in M([0,1])$ such that:
\begin{equation}\label{inde}
\langle Q \gl,\mu\rangle \, = \, 0, \qquad \langle Q\gl,\gl \rangle
+ \langle Q\mu,\mu\rangle\, = \, 1.
\end{equation}
We set for all $\go\in C([0,1])$:
\[
\gga(\go) \, := \, \int_0^1 \go_s \, \gl(ds), \qquad \Lambda_t \, :=
\, Q\gl(t), \ t\in[0,1], \qquad I \, := \, \langle Q\gl,\gl \rangle,
\]
\[
a(\go) \, := \, \int_0^1 \go_s \, \mu(ds), \qquad M_t \, := \,
Q\mu(t), \ t\in[0,1], \qquad 1-I \, = \, \langle Q\mu,\mu \rangle,
\]
and we notice that $\gga(X)\sim N(0,I)$, $a(X)\sim N(0,1-I)$ and
$\{\gga(X),a(X)\}$ are independent by \eqref{inde}. We fix a
constant $\kappa\in \bbR$ and if $I<1$ we define the continuous
processes
\[
Y_t\, := \, X_t+(\Lambda_t+M_t)\left(\kappa-a(X)-\gga(X) \right),
\quad t\in[0,1],
\]
\[
Z_t\, := \, X_t+\frac 1{1-I}\, M_t\left(\kappa-a(X)-\gga(X) \right),
\quad t\in[0,1].
\]
\begin{lemma}\label{abco} Suppose that $I<1$. Then
for all bounded Borel $\Phi:C([0,1])\mapsto\bbR$:
\begin{equation}\label{stga}
\bbE\left[\Phi(Y)\right] \, = \, \bbE\left[\Phi(Z) \, \rho(Z)
\right],
\end{equation}
where for all $\go\in C([0,1])$:
\[
\rho(\go) \, := \, \frac 1{\sqrt{1-I}} \, \exp\left(-\frac 12
\,\frac{1}{1-I} \, \left(\gga(\go)-\kappa\right)^2+\, \frac12 \,
\kappa^2 \right).
\]
\end{lemma}
\noindent We postpone the proof of Lemma \ref{abco} to section
\ref{secproof}.

\subsection{Proof of Theorem \ref{4.1}}

\noindent If $(X,Y)$ is a centered Gaussian vector and $Y\mapsto\R$
is not constant, then it is well known that a regular conditional
distribution of $X$ given $Y=y\in\R$ is given by the law of
\[
X - \frac{\sigma_{XY}}{\sigma_{YY}}(Y-y), \qquad \text{where} \quad
\sigma_{XY}=\E(XY), \quad \sigma_{YY} = \E(Y^2).
\]
We apply this property to $X=(\beta_t, t\in[0,1])$ and to
$Y=\int_0^1 \beta$. Notice that for all $t\in[0,1]$:
\[
\bbE\left[\beta_t \int_0^1 \beta_r \, dr\right] = \frac{t(1-t)}2,
\qquad \bbE\left[\left(\int_0^1 \beta_r \, dr\right)^2\right] =
\frac 1{12}.
\]
Therefore, for all $c\in\bbR$, a regular conditional distribution of
the law of $\beta$ conditioned to $\int_0^1 \beta=c$ is given by the
law of the process:
\begin{equation}\label{betac}
\beta^c_t \, := \, \beta_t \, + \, 6\, t\, (1-t) \left(c-\int_0^1
\beta\right), \qquad t\in[0,1].
\end{equation}
\begin{lemma}\label{l44}
Let $c\in\bbR$. For all bounded Borel $\Phi:C([0,1])\mapsto\bbR$:
\[
\bbE\left[\Phi(\beta) \, \left| \, \int_0^1\beta=c \right. \right]
\, = \, \bbE\left[\Phi(\beta^c)\right] \, = \,
\bbE\left[\Phi\left(\Gamma^\beta\right)\,
\rho_1\left(\Gamma^\beta\right)\right]
\]
where for all $\go\in C([0,1])$
\begin{equation}\label{gggg}
\Gamma^\go_t \, = \, \left\{ \begin{array}{ll} {\displaystyle
\go_t, \qquad \qquad \qquad \qquad \qquad \qquad \qquad  \quad
t\in[0,1/3]\cup[2/3,1] }
\\ \\
{\displaystyle \go_t+ 18\left(9\, t\, (1-t)-2\right)\left(c-\int_0^1
\go\right), \quad t\in[1/3,2/3]}
\end{array} \right.
\end{equation}
\[
\rho_1(\go) \,:= \, \sqrt{27} \, \exp\left( -162\left(\int_0^\frac
13 \left(\go_r+\go_{1-r}\right) dr + \frac{\go_\frac 13 +\go_\frac
23}6- c\right)^2 +6\, c^2 \right).
\]
\end{lemma}
\noindent {\bf Proof}. We shall show that we are in the situation of
Lemma \ref{abco} with $X=\beta$, $Y=\beta^c$ and $Z=\Gamma^\beta$.
In the notation of Lemma \ref{abco}, we consider
\[
\gl(dt) \, := \, \sqrt{12} \left( 1_{[0,\frac 13]\cup[\frac
23,1]}(t) \, dt + \frac{\delta_{\frac 13}(dt)+\delta_{\frac
23}(dt)}6 \right),
\]
\[
\mu(dt) \, := \, \sqrt{12} \left( 1_{[\frac 13,\frac 23]}(t) \, dt -
\frac{\delta_{\frac 13}(dt)+\delta_{\frac 23}(dt)}6 \right),
\]
and $\kappa:=\sqrt{12} \, c$. Then:
\[
\gga(\beta) \, = \, \sqrt{12}\int_0^{\frac 13} \left(\beta_r+\frac
12 \, \beta_{\frac 13}\right) dr + \sqrt{12}\int_{\frac 23}^1
\left(\beta_r+\frac 12 \, \beta_{\frac 23}\right) dr, \quad I \, =
\, \frac{26}{27}
\]
\[
a(\beta) \, = \, \sqrt{12}\int_{\frac 13}^{\frac 23}
\left(\beta_r-\frac 12 \, \beta_{\frac 13} - \frac 12 \,
\beta_{\frac 23}\right) dr,
\]
\[
\Lambda_t \, = \ 1_{[0,\frac 13]\cup[\frac 23,1]}(t) \ \sqrt3 \,
t(1-t) \ + 1_{(\frac 13,\frac 23)}(t) \ \frac{2\sqrt{3}}9, \quad M_t
\, = \ 1_{[\frac 13,\frac 23]}(t) \ \sqrt3 \, t(1-t).
\]
The desired result follows by tedious direct computations and from
Lemma \ref{abco}. \qed

\begin{lemma}\label{l3.5} For all bounded Borel
$\Phi:C([0,1])\mapsto\bbR$ and $f:\bbR\mapsto\bbR$
\begin{equation}\label{e4.1}
\bbE\left[\Phi(e) \, f(\langle e,1\rangle)\right] \, = \,
\int_0^\infty 27\sqrt{\frac6{\pi^3}} \,
\bbE\left[\Phi\left(V^c\right)\, \rho^c(V^c) \, \,
1_{K_0}(V^c)\right] \, f(c) \, dc
\end{equation}
\end{lemma}
\medskip\noindent
{\bf Proof}. Define $\{B,b,\hat B\}$, processes such that:
\begin{enumerate}
\item $B$ and $\hat B$ are independent copies of a standard Brownian motion over
$[0,1/3]$
\item conditionally on $\{B,\hat B\}$, $b$ is a Brownian bridge over $[1/3,2/3]$ from
$B_{1/3}$ to $\hat B_{1/3}$.
\end{enumerate}
We set:
\[
r_t \, := \, \left\{ \begin{array}{ll} B_t \qquad \qquad t\in[0,1/3]
\\ \\
b_t \qquad \qquad t\in[1/3,2/3]
\\ \\
\hat B_{1-t} \qquad \qquad t\in[2/3,1].
\end{array} \right.
\]
Moreover we set, denoting the density of $N(0,t)(dy)$ by $p_t(y)$:
\[
\rho_2(\go) \, := \, \frac{p_\frac 13(\go_\frac  23-\go_\frac
13)}{p_1(0)} \, = \,  \sqrt{3} \, \exp\left( - \frac 32 \,
(\go_{\frac 23}-\go_{\frac 13})^2\right), \qquad \go\in C([0,1]).
\]
By the Markov property of $\beta$:
\[
\bbE \left[\Phi(r) \, \rho_2(r)\right] \, = \, \bbE[\Phi(\beta)].
\]
Then, recalling the definition of $\rho^c$ above, by Lemma
\ref{abco} and Lemma \ref{l44}:
\[
\bbE \left[\Phi(\beta^c)\right] \, = \,
\bbE\left[\Phi\left(\Gamma^\beta\right)\,
\rho_1\left(\Gamma^\beta\right)\right] \, = \,
\bbE[\Phi\left(\Gamma^r\right)\, \rho_1\left(\Gamma^r\right) \,
\rho_2(\Gamma^r)] \, = \, 9 \ \bbE[\Phi(\Gamma^r) \,
\rho^c(\Gamma^r)] \, e^{6c^2}.
\]
We recall now that $\bbP(\beta\in K_\gep)=1-\exp(-2\, \gep^2) \sim
2\, \gep^2$ as $\gep\to0$, where $K_\gep=\{\go\in C([0,1]): \go\geq
-\gep\}$. We want to compute the limit of $\frac 1{2\, \gep^2} \,
\bbE\left[\Phi(\beta^c) \, 1_{K_\gep}(\beta^c)\right]$ as $\gep\to
0$. On the other hand $\bbP(B_t\geq -\gep, \forall t\in[0,1/3])\sim
\sqrt{\frac6\pi}\, \gep$ by (\ref{reflection}).  Then by
(\ref{scaling}) and (\ref{reflection})
\begin{equation}\label{ka}
 \frac 1{2\, \gep^2} \,
\bbE\left[\Phi(\beta^c) \, 1_{K_\gep}(\beta^c)\right] \, \to \,
\frac{27}\pi \, \bbE\left[\Phi\left(V^c\right)\, \rho^c(V^c) \, \,
1_{K_0}(V^c)\right] e^{6c^2}.
\end{equation}
On the other hand, $\beta$ conditioned on $K_\gep$ tends in law to
the normalized Brownian excursion $(e_t,t\in[0,1])$, as proven in
\cite{dim}. Then we have for all bounded $f\in C(\bbR)$:
\[
\frac 1{2\, \gep^2} \, \bbE\left[\Phi(\beta) \, 1_{K_\gep}(\beta) \,
f(\langle \beta,1\rangle)\right] \, \to \, \bbE\left[\Phi(e) \,
f(\langle e,1\rangle)\right]
\]
Comparing the two formulae for all $f\in C(\bbR)$ with compact
support:
\begin{align*}
& \frac 1{2\, \gep^2} \, \bbE\left[\Phi(\beta) \, 1_{K_\gep}(\beta)
\, f(\langle \beta,1\rangle)\right] \, = \, \int_\bbR \frac 1{2\,
\gep^2} \, \bbE\left[\Phi(\beta^c) \, 1_{K_\gep}(\beta^c) \right]
f(c) \, N(0,1/12)(dc)
\\ & \to \, \int_0^\infty 27\sqrt{\frac6{\pi^3}} \, \bbE\left[\Phi\left(V^c\right)\, \rho^c(V^c) \,
\, 1_{K_0}(V^c)\right] \, f(c) \, dc \, = \, \bbE\left[\Phi(e) \,
f(\langle e,1\rangle)\right]
\end{align*}
and (\ref{e4.1}) is proven. \qed

\medskip\noindent {\it Proof of Theorem \ref{4.1}}. It only remains
to prove the positivity assertion about the density. Notice that
a.s. $V^c_t\geq 0$ for all $t\in[0,1/3]\cup[2/3,1]$, since a.s.
$m\geq 0$: therefore a.s.
\[
\{V^c_t\geq 0, \ \forall t\in[0,1]\} \, = \, \{V^c_t\geq 0, \
\forall t\in[1/3,2/3]\}.
\]
The probability of this event is positive for all $c>0$ while it is
0 for $c=0$, since $\int_0^1 V^0_t \, dt = 0$. In particular
$p_{\langle e,1\rangle}(0)=0$. Finally, $p_{\langle
e,1\rangle}(c)>0$ yields also $\cZ_c>0$ if $c>0$. The other results
follow from Lemma \ref{l3.5}. \qed

\section{The linear equation}

We start with the linear Cahn-Hilliard equation, written in abstract form:
\begin{equation}\label{OU}
\left\{ \begin{array}{ll} dZ_t = - \, A \, A_D \, Z \, dt +  B \,
dW_t,
\\ \\ Z_0(x) \, = \, x\in L^2(0,1),
\end{array} \right.
\end{equation}
where $W$ is a cylindrical white noise in $L^2(0,1)$ and
\[
D(B) \, := \, H^1_0(0,1), \quad B \, := \, \sqrt 2 \,
\frac{d}{d\theta}, \qquad D(B^*) \, := \, H^1(0,1), \quad B^* \, :=
\, -\sqrt 2 \,\frac{d}{d\theta},
\]
and we notice that $BB^*=-2A$. We define the strongly continuous
contraction semigroups in $L^2(0,1)$:
\begin{equation}\label{S}
S_t \, := \, e^{-tA A_D}, \qquad  S_t^* \, := \, e^{-tA_D A},
\qquad t\geq 0.
\end{equation}
We stress that $S$ and $S^*$ are dual to each other with respect to $\langle\cdot,\cdot\rangle_L$
but not necessarily with respect to $(\cdot,\cdot)_H$.
Then it is well known that $Z$ is equal to:
\[
Z_t(x) \, = \, S_tx \, + \, \int_0^t S_{t-s}
\, B \, dW_s
\]
and that this process belongs to $C([0,\infty);L^2(0,1))$. Notice that
\begin{equation}\label{constantZ}
\langle Z_t(x),1\rangle_L \, = \, \langle x,S^*_t1\rangle_L +
\int_0^t\langle B^*S^*_{t-s}1,dW_s\rangle_L \, = \, \langle
x,1\rangle_L,
\end{equation}
since $S^*_t1=1$ and $B^*S^*_t1=B^*1=0$.
In particular, the average of $Z$ is constant.
Now, the $L^2(0,1)$-valued r.v. $Z_t(x)$ has law:
\[
Z_t(x) \sim \cN\left(S_tx, Q_t\right), \qquad Q_t \, := \, \int_0^t
S_s BB^* S^*_s \, ds.
\]
Notice that:
\[
\frac d{ds} \, S_s (-A_D)^{-1} S^*_s \, =  \, S_s (2A)\,  S^*_S \, =
\, - \, S_s BB^* S^*_s \, = \, - \, \frac d{ds} \, Q_s,
\]
so that, recalling that $Q_D:=(-A_D)^{-1}$:
\[
Q_t \, = \, Q_D \, - \, S_t \, Q_D \, S^*_t, \qquad t\geq 0.
\]
If we let $t\to\infty$, then
$\langle S_th,k\rangle_L\to \langle 1,h\rangle_L \, \langle 1,k\rangle_L$ for all $h\in L^2(0,1)$.
Therefore the law of $Z_t(x)$ converges to the Gaussian
measure on $L^2(0,1)$:
\begin{equation}\label{mu}
\mu_c \, := \, \cN( c \cdot {\bf a},Q_\infty), \qquad Q_\infty \, =
\, Q_D - \frac 1{\langle Q_D1,1\rangle_L} \, Q_D1\otimes Q_D1,
\end{equation}
with covariance operator $Q_\infty$ and mean $c \cdot {\bf a} \in
L^2(0,1)$, where
\[
c = \overline x = \langle x,1 \rangle_L, \qquad {\bf a}_\theta := 6
\, \theta(1-\theta), \quad \theta\in[0,1].
\]
Notice that the kernel of $Q_\infty$ is $\{t1: \, t\in\bbR\}$ and
therefore $\mu_c$ is concentrated on the affine space $L^2_c$.
Finally, we introduce the Gaussian measure on $L^2(0,1)$:
\begin{equation}\label{overmu}
\mu \, := \, \cN(0,Q_D),
\end{equation}
recall (\ref{Q_D}). In this case, the kernel of $Q_D$ in $L^2(0,1)$
is the null space, so the support of $\mu$ is the full space
$L^2(0,1)$. The next result gives a description of $\mu$ and $\mu_c$
as laws of stochastic processes related to the Brownian bridge
$(\beta_\theta, \theta\in[0,1])$.
\begin{lemma}\label{Y}
Let $(\beta_\theta)_{\theta\in[0,1]}$ a Brownian bridge from $0$ to
$0$. Then $\mu$ is the law of $\beta$ and $\mu_c$ is the law of the
process $\beta^c$ defined in \eqref{betac}, i.e. of $\beta$
conditioned to $\{\int_0^1\beta=c\}$, $c\in\bbR$.
\end{lemma}
\noindent{\it Proof}. By \eqref{Q_D} we have that the $Q_D$ is given
by a symmetric kernel $(\theta\wedge \sigma-\theta\sigma, \
\sigma,\theta\in[0,1])$. Since $\bbE(\beta_t \beta_s)=t\wedge s-ts$,
for all $t,s\in[0,1]$, then it is well known that $\mu = \cN(0,Q_D)$
coincides with the law of $\beta$. Analogously, the covariance of
$\beta^0$ is by \eqref{betac}
\[
\E(\beta_t^0 \beta_s^0) = t\wedge s-ts -3 \, t(1-t) \, s(1-s),
\qquad t,s\in[0,1].
\]
By the definition of $Q_\infty$ in \eqref{mu}, this is easily seen
to be the kernel of $Q_\infty$, so that $\mu_0=\cN(0,Q_\infty)$ is
the law of $\beta^0$.  By the definitions of $\mu_c=\cN(c{\bf a},
Q_\infty)$, $\beta^c$ and ${\bf a}$, we find that $\mu_c$ is the law
of $\beta^c=\beta^0+c{\bf a}$. \qed

\medskip\noindent
In particular, $\mu_c$ is a regular conditional
distribution of $\mu(dx)$ given $\{\overline x=c\}$, ie:
\[
\mu_c(dx) \, = \, \mu(dx \, | \, \overline x=c) \, = \,
\mu(dx \, | \, L^2_c).
\]
Recall (\ref{gradients}). Then we have the following result:
\begin{prop}\label{dirgau}
Let $c\in\bbR$. The bilinear form:
\[
\Lambda^c(\varphi,\psi) \, :=  \, \int_H (\nabla_{H_0} \varphi ,
\nabla_{H_0} \psi )_H \, d\mu_c \, = \, \int_H \langle-A\nabla
\varphi,\nabla\psi\rangle_L \, d\mu_c, \quad \forall \
\varphi,\psi\in {\rm Exp}_A(H),
\]
is closable in $L^2(\mu_c)$ and the process $(Z_t(x):t\geq 0, x\in H_c)$ is associated with the
resulting symmetric Dirichlet form $(\Lambda^c,D(\Lambda^c))$. Moreover,
${\rm Lip}(H_c)\subset D(\Lambda^c)$
and $\Lambda^c(\varphi,\varphi)\leq [\varphi]_{{\rm Lip}(H_c)}^2$.
\end{prop}
\noindent{\bf Proof}. The proof is standard, since the process $Z$
is Gaussian: see \cite[\S 10.2]{dpz3}. However we include some
details since the interplay between the Hilbert structures of $H$
and $L^2(0,1)$ and the different role of the operators $A$ and $A_D$
can produce some confusion. The starting point is the following
integration by parts formula for $\mu$:
\begin{equation}\label{partsmu}
\int \partial_h \varphi \ d\mu \, = \, \int \langle-A_Dh,x\rangle_L \, \varphi(x) \, \mu(dx)
\end{equation}
for all $\varphi\in C^1_b(H)$ and $h\in D(A_D)$. By conditioning on
$\{\overline x=c\}$, (\ref{partsmu}) implies:
\begin{equation}\label{partsmuc}
\int \partial_{(h-\overline h)} \varphi \ d\mu_c \, = \, \int
\langle -A_Dh,x\rangle_L \, \varphi(x) \, \mu_c(dx).
\end{equation}
Let now $\varphi(x):=\exp(i\langle x,h\rangle_L)$ and
$\psi(x):=\exp(i\langle x,k\rangle_L)$, $x\in H$, $h,k\in D(A_D A)$. Then:
\[
\bbE \, [\varphi(Z_t(x))] \, = \, \exp\left( i\langle S_t^*h,x\rangle_L -
\frac 12\langle Q_th,h\rangle_L\right)
\]
and computing the time derivative at $t=0$ we obtain the generator of $Z$:
\begin{equation}\label{L}
L\varphi(x) \, = \, \varphi(x)\left[ -i\langle A_D Ah,x\rangle_L +
\langle Ah,h\rangle_L \right].
\end{equation}
Now we compute the scalar product in $L^2(\mu_c;\bbC)$ between $L\varphi$ and $\psi$:
\begin{align*}
\int L\varphi \ \overline{\psi} \, d\mu_c & = \int \left[ -i\langle
A_D Ah,x\rangle_L + \langle Ah,h\rangle_L\right] \exp(i\langle
h-k,x\rangle_L) \, \mu_c(dx)
\\ & = \int \left[-\langle Ah,h-k\rangle_L+ \langle Ah,h\rangle_L \right]
\exp(i\langle h-k,x\rangle_L) \, \mu_c(dx)
\\ & = \int \langle Ah,k\rangle_L \,
\exp(i\langle h-k,x\rangle_L) \, \mu_c(dx) = \int \langle A\nabla
\varphi , \nabla\overline\psi\rangle_L \, d\mu_c
\end{align*}
where $\overline \psi$ is the complex conjugate of $\psi$ and in the second equality we
have used (\ref{partsmuc}). It follows that $(L,{\rm Exp}_A(H))$ is symmetric in $L^2(\mu_c)$ and
the rest of the proof is standard. \quad $\square$

\section{The approximating equation}\label{ae}

We consider now the following approximating equation:
\begin{equation}\label{ap}
\left\{ \begin{array}{ll} {\displaystyle \frac{\partial
u^{\gep,\alpha}}{\partial t}=- \frac{\partial^2 }{\partial \theta^2}
\left(\frac{\partial^2 u^{\gep,\alpha}}{\partial \theta^2} +
\frac{(u^{\gep,\alpha}+\alpha)^-}\gep \right) +
\frac{\partial}{\partial\theta} \dot{W}, }
\\ \\
{\displaystyle u^{\gep,\alpha}(t,0)=u^{\gep,\alpha}(t,1)=
\frac{\partial^3 u^{\gep,\alpha}}{\partial
\theta^3}(t,0)=\frac{\partial^3 u^{\gep,\alpha}}{\partial
\theta^3}(t,1)=0 }
\\ \\
u^{\gep,\alpha}(0,\theta)=x(\theta)
\end{array} \right.
\end{equation}
where $\gep>0$. Notice that this is a monotone gradient system in
$H$: see \cite[Chapter 12]{dpz3}, i.e. \eqref{ap} can be written as
follows,
\[
dX^{\gep,\alpha}_t = - \, A \, \left(A_D \, X^{\gep,\alpha} -\nabla
U_{\gep,\alpha}(X^{\gep,\alpha})\right)\, dt + B \, dW_t, \qquad
X_0^{\gep,\alpha}(x) \, = \, x,
\]
where $U_{\gep,\alpha}:H\mapsto \bbR^+$ is defined by
\[
U_{\gep,\alpha}(x) := \, \left\{ \begin{array}{ll} {\displaystyle
\frac{\|(x+\alpha)^-\|^2_L}\gep}, \quad &\text{if } x\in L^2(0,1)
\\ \\ +\infty, & \text{otherwise}.
\end{array} \right.
\]
We define the probability measure on $L^2(0,1)$:
\[
\nu_c^{\gep,\alpha}(dx) \, := \, \frac 1{Z^{\gep,\alpha}_c} \,
\exp\left(-\, U_{\gep,\alpha}(x) \right)\ \mu_c(dx),
\]
where $Z^{\gep,\alpha}_c$ is a normalization constant. Now,
recalling (\ref{gradients}), we introduce the symmetric bilinear
form:
\[
{\cE}^{\varepsilon,\alpha,c}(\varphi,\psi) \, := \, \int_H
(\nabla_{H_0} \varphi , \nabla_{H_0} \psi )_H \,
d\nu^{\gep,\alpha}_c \, = \, \int_H \langle-A\nabla
\varphi,\nabla\psi\rangle_L \, d\nu^{\gep,\alpha}_c, \quad \forall \
\varphi,\psi\in {\rm Exp}_A(H).
\]
Notice that this symmetric form is naturally associated with the
operator:
\begin{equation}\label{L^gep}
L^{\gep,\alpha} \varphi(x) \, := \, L\varphi(x) \, + \, \langle
\nabla U_{\gep,\alpha}(x),A\nabla\varphi\rangle_L, \qquad \forall \
\varphi\in {\rm Exp}_A(H), \ x\in L^2(0,1),
\end{equation}
where $L\varphi$ is defined in (\ref{L}) above.
The following proposition states that equation (\ref{ap}) has a unique martingale
solution, associated with the Dirichlet form arising from the closure of $({\cE}^{\varepsilon,c},{\rm Exp}_A(H))$.
Moreover, it states that the associated semigroup is Strong Feller.
\begin{prop}\label{apprteo} Let $c\in\bbR$ and $\gep> 0$.
\begin{enumerate}
\item $(L^{\gep,\alpha},{\rm Exp}_A(H))$ is essentially
self-adjoint in $L^2(\nu^{\gep,\alpha}_c)$
\item $({\cE}^{\varepsilon,\alpha,c},{\rm Exp}_A(H))$ is
closable in $L^2(\nu^{\gep,\alpha}_c)$: we denote by
$(\cE^{\gep,\alpha,c},D(\cE^{\gep,\alpha,c}))$ the closure. Moreover
${\rm Lip}(H_c)\subset D(\cE^{\gep,\alpha,c})$ and
$\cE^{\gep,\alpha,c}(\varphi,\varphi)\leq [\varphi]_{{\rm
Lip}(H_c)}^2$.
\end{enumerate}
\end{prop}
\noindent For the proof, see \cite{dpz3} and \S 9 of \cite{dpr}.
\qed

\section{Convergence of the stationary measures}

The first technical result is the convergence of
$\nu^{\gep,\alpha}_c$ as $\gep\to0^+$ and then $\alpha\to0^+$, and
in particular the tightness in a suitable H\"older space. By Lemma
\ref{Y}, $\mu_c$ is the law of $\beta^c$ defined in \eqref{betac}.
We set $K_\alpha=\{\go\in C([0,1]): \go\geq -\alpha\}$ and for
$\alpha>0$
\[
\nu^{0,\alpha}_c := \mu_c( \, \cdot \, | \, K_\alpha) = \text{law of
} \beta^c \text{ condiditioned to be greater or equal to } -\alpha.
\]
This is well defined, since $\mu_c(K_\alpha)>0$, and it is easy to
see that t
\begin{equation}\label{tig0}
\nu^{\gep,\alpha}_c \to \nu^{0,\alpha}_c \quad \text{as } \gep\to 0,
\qquad \text{weakly in } C([0,1]).
\end{equation}
Moreover, since $\beta^c$ has the same path regularity as $\beta$,
it is easy to see that for all $\alpha>0$, $\gamma\in(0,1/2)$ and
$r\geq 1$:
\begin{equation}\label{tig1}
\sup_{\gep>0} \left( \int_H \| x\|^p_{W^{\gamma,r}(0,1)} \,
d\nu_c^{\gep,\alpha}(x) \right)^{\frac1p} \, \leq \, c_{\alpha}
\left( \int_H \| x\|^p_{W^{\gamma,r}(0,1)} \, d\mu_c(x)
\right)^{\frac1p} < + \infty.
\end{equation}
We also need a similar tightness and convergence result for
$(\nu^{0,\alpha}_c)_{\alpha>0}$. We recall the definition
\[
\nu_c := \bbP[ e\in \, \cdot \, | \langle e,1\rangle = c ], \quad
\text{as defined in Theorem \ref{4.1}}.
\]
\begin{lemma}\label{coti}
As $\alpha\to0^+$, $\nu^{0,\alpha}_c$ coverges weakly in $C([0,1])$
to $\nu_c$ and
\begin{equation}\label{tig2}
\sup_{\alpha>0} \left( \int_H \| x\|^p_{W^{\gamma,r}(0,1)} \,
d\nu_c^{0,\alpha}(x) \right)^{\frac1p} \,  < + \infty.
\end{equation}
\end{lemma}
\noindent{\it Proof}. We use Lemma \ref{l44}. We recall that $\bbE
\left[\Phi(\beta^c)\right] \, = \,
\bbE\left[\Phi\left(\Gamma^\beta\right)\,
\rho_1\left(\Gamma^\beta\right)\right] $ for all bounded Borel
$\Phi:C([0,1])\mapsto\bbR$. Moreover, as proven in the Proof of
Lemma \ref{l3.5}, the law of $\Gamma^\beta$ conditioned on
$K_\alpha$ converges to the law of $V^c$, defined in \eqref{v2},
i.e.
\[
\lim_{\alpha\to0^+} \bbE\left[ \Phi( \Gamma^\beta ) \, | \,
\Gamma^\beta \in K_\alpha\right] = \bbE[ \Phi(V^c) ].
\]
Notice that $\rho_1$ is positive, continuous on $C([0,1])$ and
bounded by a constant. Then we have
\begin{align*}
& \bbE\left[ \|\beta^c \|^p_{W^{\gamma,r}(0,1)} \, | \ \beta^c\in
K_\alpha \right] = \frac{\bbE\left[ \| \Gamma^\beta
\|^p_{W^{\gamma,r}(0,1)} \, \rho_1(\Gamma^\beta) \, 1_{(\Gamma^\beta
\in K_\alpha)}\right]}{\bbE\left[ \rho_1(\Gamma^\beta) \,
1_{(\Gamma^\beta \in K_\alpha)}\right]}
\\ & \leq \kappa_1 \, \bbE\left[ \| \Gamma^\beta
\|^p_{W^{\gamma,r}(0,1)} \, | \, \Gamma^\beta \in K_\alpha\right]
\cdot \frac1{\bbE\left[ \rho_1(\Gamma^\beta) \, | \, \Gamma^\beta
\in K_\alpha \right]} \leq \kappa_2 \bbE\left[ \| \Gamma^\beta
\|^p_{W^{\gamma,r}(0,1)} \, | \, \Gamma^\beta \in K_\alpha\right],
\end{align*}
where the last inequality follows from the convergence $\bbE\left[
\rho_1(\Gamma^\beta) \, | \, \Gamma^\beta \in K_\alpha \right] \to
\bbE\left[ \rho_1(V^c)  \right]>0$, $\alpha\to0^+$. Then it only
remains to prove that
\begin{equation}\label{dennum}
\sup_{\alpha>0} \frac{\bbE\left[ \| \Gamma^\beta
\|^p_{W^{\gamma,r}(0,1)} \, 1_{(\Gamma^\beta \in K_\alpha)}
\right]}{\bbP\left[ \Gamma^\beta \in K_\alpha\right]} <+ \infty.
\end{equation}

We start with the numerator. We fix three functions
$\phi_i:[0,1]\mapsto \bbR_+$ of class $C^\infty$, such that
$\phi_1+\phi_2+\phi_3\equiv 1$, the support of $\phi_1$ is in
$[0,1/3)$, the support of $\phi_3$ is in $(2/3,1]$ and the support
of $\phi_2$ is in $(1/6,5/6)$. Then it is enough to estimate
\[
\bbE\left[ \| \varphi_i \cdot \Gamma^\beta \|^p_{W^{\gamma,r}(0,1)}
\, 1_{(\Gamma^\beta \in K_\alpha)}\right], \qquad i=1,2,3.
\]
Notice that $\varphi_1\Gamma^\beta=\varphi_1\beta$. We set
$I=[0,1/3]$ and we denote by $(\beta^{0,a}_\theta, \theta\in I)$,
resp. $(m^{b,a}_\theta, \theta\in I)$, the Brownian bridge from $0$
to $a$ over the interval $I$, respectively the 3-dimensional Bessel
bridge from $b$ to $a$ over the interval $I$. Then, denoting by
$p_t$ the density of $\cN(0,t)$,
\begin{align*}
& \E(\Phi(\beta_\theta, \theta\in I) \, | \, \beta\geq -\alpha
\text{ on } I ) = \int_{-\alpha}^\infty \E(\Phi(\beta_\theta^{0,a},
\theta\in I) \, | \, \beta^{0,a}\geq -\alpha \text{ on } I ) \,
p_{2/9}(a) \, da
\\ & = \int_{-\alpha}^\infty \E(\Phi(m_\theta^{\alpha,a+\alpha}-\alpha,
\theta\in I) ) \, p_{2/9}(a) \, da
\end{align*}
where in the former equality we use the Markov property of $\beta$
and in latter the equality in law between Brownian bridges
conditioned to be positive and  3-dimensional Bessel bridges. Then
\begin{align*}
& \bbE\left[ \| \varphi_1 \cdot \Gamma^\beta
\|^p_{W^{\gamma,r}(0,1)} \, 1_{(\Gamma^\beta \in K_\alpha)}\right]
\leq \bbE\left[ \| \varphi_1 \cdot \beta \|^p_{W^{\gamma,r}(0,1)} \,
1_{(\beta\geq -\alpha \text{ on }[0,1/3]\cup[2/3,1])}\right]
\\ & = \int_{-\alpha}^\infty \bbE\left[ \| \varphi_1 \cdot
\left(m^{\alpha,a+\alpha}-\alpha\right) \|^p_{W^{\gamma,r}(0,1)}
\right] \, p_{2/9}(a)\cdot \gamma_\alpha(a) \, da,
\end{align*}
where $\gamma_\alpha(a)=\bbP(\beta\geq -\alpha \text{ on
}[0,1/3]\cup[2/3,1] \ | \ \beta_{1/3}=a)$. Then it is easy to
conclude that
\[
\sup_{\alpha>0} \frac1{\alpha^2} \, \bbE\left[ \| \varphi_1 \cdot
\Gamma^\beta \|^p_{W^{\gamma,r}(0,1)} \, 1_{(\Gamma^\beta \in
K_\alpha)}\right] < + \infty.
\]
By symmetry, the same estimate holds for $\varphi_3 \cdot
\Gamma^\beta$. As for $\varphi_2 \cdot \Gamma^\beta$, conditioning
on the values of $\beta_{1/3}$ and $\beta_{2/3}$ and using an
analogous argument, we find similarly that
\[
\sup_{\alpha>0} \frac1{\alpha^2} \, \bbE\left[ \| \varphi_2 \cdot
\Gamma^\beta \|^p_{W^{\gamma,r}(0,1)} \, 1_{(\Gamma^\beta \in
K_\alpha)}\right] < + \infty.
\]

We estimate now the denominator of the r.h.s. of \eqref{dennum}.
Recall the definition \eqref{gggg} of $\Gamma^\omega$ for $\omega\in
C([0,1])$. Notice that
\[
\int_0^1 \omega \geq c \ \Longrightarrow \ \Gamma^\omega_t \leq
\omega_t, \quad \forall \ t\in[0,1],
\]
since $9t(1-t)-2\geq 0$ for all $t\in[1/3,2/3]$. This means that
\begin{align*}
& \bbP(\Gamma^\beta \in K_\alpha) \geq \bbP\left(\Gamma^\beta \in
K_\alpha, \ \int_0^1 \beta\geq c\right) \geq \bbP\left(\beta \in
K_\alpha, \ \int_0^1 \beta\geq c\right)
\\ & = \bbP\left(\int_0^1 \beta\geq c \ \Big| \ \beta \in
K_\alpha\right) \cdot \bbP\left(\beta \in K_\alpha\right) \sim
\bbP\left(\int_0^1 e\geq c \right) \, 2\alpha^2, \quad \alpha\to0^+,
\end{align*}
since $\bbP\left(\int_0^1 \beta\geq c \ | \ \beta \in
K_\alpha\right) \to \bbP\left(\int_0^1 e\geq c \right)>0$. Then
\eqref{dennum} is proven.

In order to show that $\nu^{0,\alpha}_c$ indeed converges to
$\nu_c$, it is enough to recall formula \eqref{ka} above and the
second result of Theorem \ref{4.1}. \qed

\section{A general convergence result}\label{asz}

In this section we recall two results of \cite{asz}, which we shall
apply in section \ref{sec:weak} to the convergence in law of the
solutions of \eqref{ap} to the solution of \eqref{00}. These
processes are reversible and associated with a gradient-type
Dirichlet form. Moreover their invariant measures (respectively,
$\nu^{\gep,\alpha}_c$ and $\nu_c$, are {\it log-concave}; a
probability measure $\gamma$ on $H$ is log-concave if for all pairs
of open sets $B,\, C\subset H$
\begin{equation}\label{deflogconc}
\log\gamma\left((1-t)B\right)\geq (1-t)\log\gamma(B)+t\log\gamma(C)
\qquad\forall t\in (0,1).
\end{equation}
If $H=\R^k$, then the class of log-concave probability measures
contains all measures of the form (here $\cL_k$ stands for Lebesgue
measure)
\begin{equation}
  \label{eq:basic_example}
  \gamma :=  \frac 1Z \, e^{-V} \cL_k,
\end{equation}
where $V:H=\R^k\to\R$ is convex and $Z:= \int_{\R^k} e^{-V} \,
dx<+\infty$, see Theorem~9.4.11 in \cite{ags}, in particular all
Gaussian measures. Notice that the class of log-concave measures is
closed under weak convergence. Therefore, it is easy to see by an
approximation argument that $\nu^{\gep,\alpha}_c$ and $\nu_c$ are
log-concave.

We denote by $X_t:H^{[0,+\infty[}\to H$ the coordinate process
$X_t(\omega):=\omega_t$, $t\geq 0$. Then we recall one of the main
results of \cite{asz}. We notice that the support of $\nu_c$ in $H$
is $K_c$, the closure in $H$ of $\{h\in L^2(0,1): h\geq 0, \ \langle
h,1\rangle_L =c \}$, and the closed affine hull in $H$ of $K_c$ is
$H_c$.
\begin{prop}[Markov process and Dirichlet form associated with $\nu_c$ and
$\|\cdot\|_{H_0}$]\label{main1} $ $
\begin{itemize}
\item[(a)]
The bilinear form ${\cE}={\cE}_{\nu_c,\|\cdot\|_{H_0}}$ given by
\begin{equation}\label{diri}
{\cE}(u,v) \, := \, \int_{K_c} (\nabla_{H_0} u,\nabla_{H_0} v)_{H}
\, d\nu_c, \qquad u,\,v\in C^1_b(H_c),
\end{equation}
is closable in $L^2(\nu_c)$ and its closure $(\cE,D(\cE))$ is a
symmetric Dirichlet Form. Furthermore, the associated semigroup
$(P_t)_{t\geq 0}$ in $L^2(\nu_c)$ maps $L^\infty(\nu_c)$ in
$C_b(K_c)$.
\item[(b)] There exists a unique Markov family $(\bbP_x:x\in K_c)$ of probability
measures on $K_c^{[0,+\infty[}$ associated with $\cE$. More
precisely, $\E_x[f(X_t)]=P_tf(x)$ for all bounded Borel functions
and all $x\in K_c$. \smallskip
\item[(c)] For all $x\in K_c$, $\bbP_x^*\left(C(]0,+\infty[;H)\right)=1$ and
$\E_x[\|X_t-x\|^2]\to 0$ as $t\downarrow 0$. Moreover,
$\bbP_x^*\left(C([0,+\infty[;H)\right)=1$ for $\nu_c$-a.e. $x\in
K_c$.
\smallskip
\item[(d)] $(\bbP_x:x\in K_c)$ is reversible with respect to $\nu_c$,
i.e. the transition semigroup $(P_t)_{t\geq 0}$ is symmetric in
$L^2(\nu_c)$; moreover $\nu_c$ is invariant for $(P_t)$, i.e.
$\nu_c(P_tf)=\nu_c(f)$ for all $f\in C_b(K_c)$ and $t\geq 0$.
\end{itemize}
\end{prop}

\medskip\noindent
Let $(\bbP_x^{\gep,\alpha,c}:x\in H_c)$ (respectively $(\bbP_x:x\in
K_c)$) be the Markov process in $[0,+\infty[^{H_c}$ associated to
(resp. in $[0,+\infty[^{K_c}$ associated to $\nu_c$) given by
Proposition~\ref{main1}. We denote by $\bbP_{\nu_c^N}^N:=\int
\bbP_{x}^N \, d\nu_c^N(x)$ (resp. $\bbP_{\nu_c}:= \int \bbP_{x} \,
d\nu_c(x)$) the associated stationary measures.

With an abuse of notation, we say that a sequence of measures $({\bf
P}_n)$ on $C([a,b];H)$ converges weakly in $C([a,b];H_w)$ if, for
all $m\in\N$ and $h_1,\ldots,h_m\in H$, the process $(\langle
X_\cdot,h_i\rangle_H, \, i=1,\ldots,m)$ under $({\bf P}_n)$
converges weakly in $C([a,b];\R^m)$ as $n\to\infty$.

In this setting we have the following stability and tightness
result:
\begin{theorem}[Stability and tightness]\label{main3}
Then, for any $x\in K_c$ and $0<\varepsilon\leq T<+\infty$,
\[
\lim_{\alpha\to0^+} \lim_{\gep\to0^+}
\bbP_{x}^{\gep,\alpha,c}=\bbP_x, \quad \text{weakly in
$C([\varepsilon,T];H_w)$}.
\]
\end{theorem}
\noindent{\bf Proof}. This result follows from Theorem 1.5 in
\cite{asz}, where it is stated that the weak convergence of the
invariant measures of a sequence of processes as in Proposition
\ref{main1} implies the weak convergence of the associated
processes. Since $\lim_{\alpha\to0^+} \lim_{\gep\to0^+}
\nu^{\gep,\alpha}_c=\nu_c$, we obtain the result. \qed

\section{Existence of weak solutions of equation \eqref{00}}\label{sec:weak}

In this section we prove the following result on weak existence of
solutions to equation \eqref{00}. We define the Polish space
$E_T:=C(O_T)\times  M(O_T)\times C(O_T)$, where $O_T:=]0,T]\times
[0,1]$ and $M(O_T)$ is the space of all locally finite positive
measures on $]0,T]\times ]0,1[$, endowed with the topology of
convergence on compacts in $]0,T]\times ]0,1[$.

\begin{prop}\label{weak}
Let $c>0$, $u_0=x\in K_c$ and $u^{\gep,\alpha}$ the solution of
\eqref{ap}. Set $\eta^{\gep,\alpha}\in M(O_T)$,
\[
\eta^{\gep,\alpha}(dt,d\theta) :=
\frac{(u^{\gep,\alpha}(t,\theta)+\alpha)^-}\gep \, dt \, d\theta.
\]
Then $(u^{\gep,\alpha},\eta^{\gep,\alpha},W)$ converges in law to
$(u,\eta,W)$, stationary weak solution of (\ref{00}), in $E_T$, for
any $T\geq 0$. The law of $u$ is $\bbP_x$ and therefore $(u,
u_0=x\in K_c)$ is the Markov process associated with the Dirichlet
form \eqref{diri}.
\end{prop}

\noindent We shall use the following easy result:
\begin{lemma}\label{bouv}
Let $\zeta(dt,d\theta)$ be a finite signed measure on
$[\delta,T]\times [0,1]$ and $v\in C([\delta,T]\times [0,1])$.
Suppose that for all $s\in[\delta,T]$:
\begin{equation}\label{bouve}
\int_{[s,T]\times [0,1]} h_\theta \, \zeta(dt,d\theta) \, = \, 0,
\qquad \forall \ h\in C([0,1]), \ \overline h=0,
\end{equation}
and
\begin{equation}\label{bouve1}
\overline v_s \, = \, c>0, \qquad \int_{[s,T]\times [0,1]} v \,
d\zeta \, = \, 0.
\end{equation}
Then $\zeta\equiv 0$.
\end{lemma}
\noindent{\bf Proof}. Setting $h:=k-\overline k$, $k\in C([0,1])$,
we obtain by (\ref{bouve}) for all $\delta\leq s \leq t\leq T$:
\[
\int_0^1 k_\theta \, \zeta([s,t]\times  d\theta) \, = \,
\zeta([s,t]\times [0,1]) \int_0^1 k_\theta\, d\theta, \qquad \forall
\, k\in C([0,1]).
\]
This implies $\zeta(dt,d\theta) \, = \, \gamma(dt) \, d\theta$,
where $\gamma(t):=\zeta([\delta,t]\times [0,1])$, $t\in[\delta,T]$,
is a process with bounded variation. Then by (\ref{bouve1}):
\[
0 \, = \, \int_{[s,t]\times [0,1]} v \, d\zeta \, = \, \int_s^t
\left(\int_0^1 v_s(\theta) \, d\theta\right) \, \gamma(ds) \, = \, c
\, (\gamma(t)-\gga(s)),
\]
i.e. $\gamma(t)-\gamma(s)=0$, since $c>0$. \qed

\noindent{\bf Proof of Proposition \ref{weak}}. Recall that
$\bbP_{x}^{\gep,\alpha,c}$ is the law of $u^{\gep,\alpha}$ if
$u^{\gep,\alpha}_0=x$. By Theorem \ref{main3} and Skorohod's Theorem
we can find a probability space and a sequence of processes
$(v^\gep,w^\gep)$ such that $(v^\gep,w^\gep)\to(v,w)$ in $C(O_T)$
almost surely and $(v^\gep,w^\gep)$ has the same distribution as
$(u^{\gep},W)$ for all $\gep>0$, where $O_T:=]0,T]\times [0,1]$.
Notice that $v\geq 0$ almost surely, since for all $t$ the law of
$v_t(\cdot)$ is $\gamma$ which is concentrated on $K$ and moreover
$v$ is continuous on $O_T$. We set now:
\[
\eta^\gep(dt,d\theta) \, := \, \frac 1{\gep} \,
f\left(v^\gep_t(\theta)\right) \, dt\, d\theta.
\]
From (\ref{ap}) we obtain that a.s. for all $T\geq 0$ and $h\in
D(A^2)$ and $\overline h=0$:
\begin{equation}\label{exlim}
\exists \lim_{\gep\to 0^+} \int_{O_T} h_\theta \,
\eta^\gep(dt,d\theta).
\end{equation}
The limit is a random distribution on $O_T$. We want to prove that
in fact $\eta^\gep$ converges as a measure in the dual of $C(O_T)$
for all $T\geq 0$. For this, it is enough to prove that the mass
$\eta^\gep(O_T)$ converges as $n\to\infty$.

Suppose that $\{\eta^\gep(O_T)\}_n$ is unbounded. We define
$\zeta^\gep:=\eta^\gep/\eta^\gep(O_T)$. Then $\zeta^\gep$ is a
probability measure on the compact set $O_T$. By tightness we can
extract from any sequence $\gep_n\to 0$ a subsequence along which
$\zeta^\gep$ converges to a probability measure $\zeta$. By the
uniform convergence of $v^\gep$ we can see that the contact
condition $\int_{O_T} v \, d\zeta=0$ holds. Moreover, dividing
(\ref{ap}) by $\eta^\gep(O_T)$ for $t\in[0,T]$, we obtain that
$\int_{O_t} h_\theta \, \zeta(ds,d\theta) = 0$ for all $h\in D(A^2)$
with $\overline h=0$ and by density for all $h\in C([0,1])$ with
$\overline h=0$.

Then $\zeta$ and $v$ satisfy (\ref{bouve}) and (\ref{bouve1}) above,
and therefore by Lemma \ref{bouv}, $\zeta\equiv 0$, a contradiction
since $\zeta$ is a probability measure. Therefore
$\limsup_{n\to\infty} \eta^\gep(O_T)<\infty$.

By tightness, for any subsequence in $\bbN$ we have convergence of
$\eta^\gep$ to a finite measure $\eta$ on $[0,T]\times [0,1]$ along
some sub-subsequence. Let $\eta_i$, $i=1,2$, be two such limits and
set $\zeta:=\eta_1 -\eta_2$. By (\ref{exlim}) and by density:
\[
\int_{O_T} h_\theta \, \eta_1(dt,d\theta) \, = \, \int_{O_T}
h_\theta \, \eta_2(dt,d\theta), \qquad \forall \, h\in C([0,1]), \
\overline h=0,
\]
i.e. $\zeta$ and $v$ satisfy (\ref{bouve}) and (\ref{bouve1}) above.
By Lemma \ref{bouv}, $\zeta\equiv 0$, i.e. $\eta_1=\eta_2$.
Therefore, $\eta^\gep$ converges as $n\to\infty$ to a finite measure
$\eta$ on $]0,T]\times [0,1]$. It is now clear that the limit
$(u,\eta,W)$ satisfies \eqref{distr}.

Finally, we need to prove that the contact condition holds, i.e.
that $\int_{(0,\infty)\times [0,1]}v \, d\eta \, = \, 0$. Since
$f\geq 0$ and $f(u)>0$ for $u>0$, then $u\, f(u)\leq 0$ for all
$u\in \bbR$. Then for any continuous positive
$\varphi:(0,1)\mapsto\bbR$ with compact support
\[
0\, \geq \, \int_{[0,T]\times [0,1]} \, \varphi\,  v^\gep \,
d\eta^\gep \, \to \, \int_{[0,T]\times [0,1]} \, \varphi\,  v \,
d\eta
\]
by the uniform convergence of $v^\gep$ to $v$ and the convergence of
$\eta^\gep$ to $\eta$ on compacts. Since $v\geq 0$ and $\eta$ is a
positive measure, then $\int_{[0,T]\times [0,1]} \, v \, d\eta\leq
0$ is possible only if $\int_{[0,T]\times [0,1]} \, v \, d\eta= 0$
\qed

\section{Conditioning the Brownian meander to have a fixed time
average}\label{meander}

In this section we prove an analog of Theorem \ref{4.1} for the
standard Brownian meander $(m_t,t\in[0,1])$. We set $\langle
m,1\rangle:=\int_0^1 m_r \, dr$, average of $m$. Let $B$ a standard
Brownian motion such that $\{m,B\}$ are independent and let $c\geq
0$ be a constant. We introduce the continuous processes:
\begin{align*}
& u_t \, := \, \left\{ \begin{array}{ll} {\displaystyle \frac
1{\sqrt 2}\, m_{2t}, \qquad \quad t\in[0,1/2] }
\\ \\
{\displaystyle \frac 1{\sqrt 2}\, m_1 + B_{t-\frac 12}, \quad
t\in[1/2,1],}
\end{array} \right.
\\
& U^c_t \, := \, \left\{ \begin{array}{ll} {\displaystyle u_t,
\qquad \quad t\in[0,1/2] }
\\ \\
{\displaystyle u_t + \left(12\, t\, (2-t)-9\right)\left(c-\int_0^1
u\right), \quad t\in[1/2,1].}
\end{array} \right.
\end{align*}
Notice that $\int_0^1 U^c_t \, dt = c$.
\begin{theorem}\label{1.1}
Setting for all $c\geq 0$
\[
p_{\langle m,1\rangle}(c) \, := \, \sqrt\frac{24}\pi \,
\bbE\left[e^{-12\left(\int_0^{1/2}(U_r^c+U_{1/2}^c)\, dr-c\right)^2}
\, 1_{\{U^c_t\geq 0, \ \forall t\in[0,1]\}}\right],
\]
and for all bounded Borel $\Phi:C([0,1])\mapsto\bbR$ and $c>0$
\[
\bbE\left[\Phi(m) \, | \, \langle m,1\rangle =c \right] \, := \,
\frac 1{\cZ_c} \ \bbE\left[\Phi(U^c)\,
e^{-12\left(\int_0^{1/2}(U_r^c+U_{1/2}^c)\, dr-c\right)^2} \,
1_{\{U^c_t\geq 0, \ \forall t\in[0,1]\}}\right],
\]
where $\cZ_c>0$ is a normalization factor, we have
\begin{enumerate}
\item $p_{\langle m,1\rangle}$ is the density of $\langle m,1\rangle$,
i.e.
\[
\bbP(\langle m,1\rangle\in dc)\, = \, p_{\langle m,1\rangle}(c)
\,1_{\{c\geq 0\}} \, dc.
\]
Moreover $p_{\langle m,1\rangle}$ is continuous on $[0,\infty)$,
$p_{\langle m,1\rangle}(c)>0$ for all $c\in(0,\infty)$ and
$p_{\langle m,1\rangle}(0)=0$.
\item $(\bbP\left[m\in\cdot \, | \, \langle m,1\rangle =c \right], c>0)$
is a regular conditional distribution of $m$ given $\langle
m,1\rangle$, i.e.
\[
\bbP(m\in\cdot\, , \langle m,1\rangle\in dc ) \, = \,
\bbP\left[m\in\cdot \, | \, \langle m,1\rangle =c \right] \
p_{\langle m,1\rangle}(c) \,1_{\{c> 0\}} \, dc.
\]
\end{enumerate}
\end{theorem}

\noindent In the notation of section \ref{secabco}, we consider
$X=(B_t,t\in[0,1])$, standard Brownian motion. It is easy to see
that for all $t\in[0,1]$:
\[
\bbE\left[B_t \int_0^1 B_r \, dr\right] = \frac{t\,(2-t)}2, \qquad
\bbE\left[\left(\int_0^1 B_r \, dr\right)^2\right] = \frac 13.
\]
Therefore, it is standard that for all $c\in\bbR$, $B$ conditioned
to $\int_0^1 B=c$ is equal in law to the process:
\[
B^c_t \, := \, B_t \, + \, \frac 32\, t\, (2-t) \left(c-\int_0^1
B\right), \qquad t\in[0,1].
\]
\begin{lemma}\label{l41}
Let $c\in\bbR$. For all bounded Borel $\Phi:C([0,1])\mapsto\bbR$:
\[
\bbE\left[\Phi(B) \, \left| \, \int_0^1 B=c\right. \right] \, = \,
\bbE\left[\Phi(B^c)\right] \, = \, \bbE\left[\Phi\left(S\right)\,
\rho(S)\right],
\]
where
\[
S_t \, := \, \left\{ \begin{array}{ll} {\displaystyle B_t, \qquad
\qquad \qquad \qquad \qquad \qquad \qquad \quad t\in[0,1/2] }
\\ \\
{\displaystyle B_t+ \left(12\, t\, (2-t)-9\right)\left(c-\int_0^1
B\right), \quad t\in[1/2,1]}
\end{array} \right.
\]
\[
\rho(\go) \, := \, \sqrt 8 \, \exp\left( -12\left(\int_0^\frac 12
\left(\go_r+\go_\frac 12\right) dr - c\right)^2 +\frac 32 \, c^2
\right), \quad \go\in C([0,1]).
\]
\end{lemma}
\noindent{\bf Proof}. We are going to show that we are in the
setting of Lemma \ref{abco} with $X=B$, $Y=B^c$ and $Z=S$. We denote
the Dirac mass at $\theta$ by $\delta_\theta$. In the notation of
section \ref{secabco}, we consider:
\[
\gl(dt) \, := \, \sqrt3 \left( 1_{[0,\frac 12]}(t) \, dt + \frac 12
\, \delta_{\frac 12}(dt)\right), \qquad \mu(dt) \, := \, \sqrt3
\left( 1_{[\frac 12,1]}(t) \, dt - \frac 12 \, \delta_{\frac
12}(dt)\right),
\]
and $\kappa:=\sqrt 3 \, c$. Then:
\[
\gga(\go) \, = \, \sqrt 3\int_0^{\frac 12} \left(\go_r+\go_{\frac
12}\right) dr, \qquad a(\go) \, = \, \sqrt 3\int_{\frac 12}^1
\left(\go_r-\go_{\frac 12}\right) dr,
\]
\[
\gga(\go)+a(\go) \, = \, \sqrt 3 \int_0^1 \go_r \, dr, \qquad I \, =
\, 3\int_0^{\frac 12} (1-r)^2 \, dr \, = \, \frac78.
\]
\[
\Lambda_t \, = \, \left\{ \begin{array}{ll} {\displaystyle \sqrt3 \
t \left(1-\frac {t}2\right), \quad t\in[0,1/2] }
\\ \\
{\displaystyle \frac{3\sqrt 3}8, \qquad \qquad t\in[1/2,1].}
\end{array} \right.
\ M_t \, = \, \left\{ \begin{array}{ll} {\displaystyle 0, \
t\in[0,1/2] }
\\ \\
{\displaystyle \sqrt3 \ t \left(1-\frac {t}2\right)-\frac{3\sqrt
3}8, \quad t\in[1/2,1].}
\end{array} \right.
\]
Tedious but straightforward computations show that with these
definitions we have $X=B$, $Y=B^c$ and $Z=S$ in the notation of
Lemma \ref{abco} and (\ref{inde}) holds true. Then the Lemma
\ref{l41} follows from Lemma \ref{abco}. \qed

\begin{lemma}\label{pro}
For all bounded Borel $\Phi:C([0,1])\mapsto\bbR$ and
$f:\bbR\mapsto\bbR$:
\[
\bbE\left[\Phi(m) \, f(\langle m,1\rangle)\right] =
 \int_0^\infty \sqrt\frac{24}\pi \,
\bbE\left[\Phi(U^c)\, e^{-12\left(\int_0^{1/2}(U_r^c+U_{1/2}^c)\,
dr-c\right)^2}\, 1_{\{U^c_t\geq 0, \ \forall t\in[0,1]\}}\right]
f(c) \, dc.
\]
\end{lemma}
\medskip
\noindent{\bf Proof}. Recall that $m$ is equal in law to $B$
conditioned to be non-negative (see \cite{dim} and \eqref{scaling}
below). We want to condition $B$ first to be non-negative and then
to have a fixed time average. It turns out that Lemma \ref{l41}
allows to compute the resulting law by inverting the two operations:
first we condition $B$ to have a fixed average, then we use the
absolute continuity between the law of $B^c$ and the law of $S$ and
finally we condition $S$ to be non-begative.

We set $K_\gep:=\{\go\in C([0,1]): \go\geq -\gep\}$, $\gep\geq 0$.
We recall that $B$ conditioned on $K_\gep$ tends in law to $m$ as
$\gep\to 0$, more generally for all $s>0$ and bounded continuous
$\Phi:C([0,s])\mapsto\bbR$, by the Brownian scaling:
\begin{equation}\label{scaling}
\lim_{\gep\to 0} \, \bbE\left[ \Phi(B_t, \, t\in[0,s]) \, \Big| \,
B_t\geq -\gep, \, \forall \, t\in[0,s] \right] \, = \, \bbE\left[
\Phi\left( {\sqrt s}\, m_{t/s}, \, t\in[0,s] \right) \right],
\end{equation}
and this is a result of \cite{dim}. By the reflection principle, for
all $s>0$:
\begin{equation}\label{reflection}
\bbP(B_t\geq -\gep, \ \forall t\in[0,s]) \, = \, \bbP(|B_s|\leq
\gep) \, \sim \, \sqrt\frac2{\pi\, s} \, \gep, \qquad \gep\to 0.
\end{equation}
In particular for all bounded $f\in C(\bbR)$
\[
\bbE\left[\Phi(m) \, f(\langle m,1\rangle)\right] \, = \,
\lim_{\gep\to 0}\, \sqrt\frac\pi2 \, \frac 1\gep \,
\bbE\left[\Phi(B) \, 1_{K_\gep}(B) \, f(\langle B,1\rangle)\right].
\]
We want to compute the limit of $\frac 1\gep \, \bbE\left[\Phi(B^c)
\, 1_{K_\gep}(B^c)\right]$ as $\gep\to 0$. Notice that $S$, defined
in Lemma \ref{l41}, is equal to $B$ on $[0,1/2]$. Therefore, by
(\ref{scaling}) and (\ref{reflection}) with $s=1/2$:
\[
\sqrt\frac\pi2 \, \frac 1\gep \, \bbE\left[\Phi(B^c) \,
1_{K_\gep}(B^c)\right] \, \to \, \sqrt 2 \, \bbE\left[\Phi(U^c)\,
\rho(U^c) \, 1_{K_0}(U^c)\right].
\]
Comparing the last two formulae for all $f\in C(\bbR)$ with compact
support:
\begin{align*}
& \sqrt\frac\pi2 \, \frac 1\gep \, \bbE\left[\Phi(B) \,
1_{K_\gep}(B) \, f(\langle B,1\rangle)\right] \, = \, \int_\bbR
\sqrt\frac\pi2 \, \frac 1\gep \, \bbE\left[\Phi(B^c) \,
1_{K_\gep}(B^c) \right] f(c) \, N(0,1/3)(dc)
\\ & \to \, \int_0^\infty \sqrt\frac{24}\pi \,
\bbE\left[\Phi(U^c)\, e^{-12\left(\int_0^{1/2}(U_r^c+U_{1/2}^c)\,
dr-c\right)^2}\, 1_{K_0}(U^c)\right]  \, f(c) \, dc
 = \bbE\left[\Phi(m) \, f(\langle m,1\rangle)\right]
\end{align*}
and the Lemma is proven. \qed

\medskip\noindent{\bf Proof of Theorem \ref{1.1}} The results follow
from Lemma \ref{pro}, along the lines of the proof of Theorem
\ref{4.1}. \qed

\medskip\noindent
It would be now possible to repeat the results of sections and prove
existence of weak solutions to the SPDE
\begin{equation}\label{00fi}
\left\{ \begin{array}{ll} {\displaystyle \frac{\partial u}{\partial
t}=- \frac{\partial^2 }{\partial \theta^2} \left(\frac{\partial^2
u}{\partial \theta^2} + \eta \right)  + \sqrt 2 \,
\frac{\partial}{\partial\theta} \dot{W}, }
\\ \\
{\displaystyle u(t,0)=\frac{\partial u}{\partial \theta}(t,1)=
\frac{\partial^3 u}{\partial \theta^3}(t,0)=\frac{\partial^4
u}{\partial \theta^4}(t,1)=0 }
\\ \\
u(0,\theta)=x(\theta)
\end{array} \right.
\end{equation}
and that such weak solutions admit $\bbP\left[m\in\cdot \, | \,
\langle m,1\rangle =c \right]$ as invariant measures, $c>0$.

\section{Proof of Proposition \ref{abco}.}
\label{secproof}

The result follows if we show that the Laplace transforms of the two
probability measures in (\ref{stga}) are equal. Notice that $Y$ is a
Gaussian process with mean $\kappa\, (\Lambda+M)$ and covariance
function:
\[
q^Y_{t,s} \, = \,
\bbE\left[\left(Y_t-\kappa\,(\Lambda_t+M_t)\right)\left(Y_s-\kappa\,(\Lambda_s+M_s)\right)\right]
\, = \, q_{t,s} \, - \, (\Lambda_t+M_t)\,(\Lambda_s+M_s),
\]
for $t,s\in[0,1]$. Therefore, setting for all $h\in C([0,1])$:
$Q_Yh(t):=\int_0^1 q^Y_{t,s} \, h_s \, ds$, $t\in[0,1]$, the Laplace
transform of the law of $Y$ is:
\[
\bbE\, \left[e^{\langle Y,h\rangle}\right] \, = \, e^{\kappa\langle
h,\Lambda+M\rangle+ \frac12 \, \langle Q_Y h,h\rangle}.
\]
Recall now the following version of the Cameron-Martin Theorem: for
all $h\in M([0,1])$
\[
{\mathbb E}\left[ \Phi(X) \, e^{\langle X,h\rangle} \right] \, = \,
e^{\frac 12\langle Qh,h\rangle} \, {\mathbb E} [\Phi(X+Qh)].
\]
Notice that $\gga(Z)=\gga(X)$, by (\ref{inde}). Therefore
$\rho(Z)=\rho(X)$. We obtain, setting $\overline
h:=h-\frac1{1-I}\langle M,h\rangle(\gl+\mu)$:
\begin{align*}
& \bbE\left[e^{\langle Z,h\rangle} \, \rho(Z) \right] \, = \,
e^{\frac\kappa{1-I}\langle M,h\rangle} \, \bbE\left[e^{\langle
X,\overline h\rangle}\, \rho(X) \right] \, = \,
e^{\frac\kappa{1-I}\langle M,h\rangle+\frac 12\langle Q\overline
h,\overline h\rangle} \, {\mathbb E} \left[\rho\left(X+Q\overline
h\right)\right] =
\\ & = \, e^{\frac\kappa{1-I}\langle M,h\rangle+\frac 12\langle Q\overline h,\overline h\rangle} \,
\frac 1{\sqrt{1-I}} \, \bbE\left[ e^{-\frac 12 \, \frac 1{1-I} \,
\left(\gga(X)+\langle \overline h,\Lambda\rangle-\kappa\right)^2+
\frac12 \, \kappa^2}\right].
\end{align*}
By the following standard Gaussian formula for $\ga\sim
N(0,\sigma^2)$, $\sigma\geq 0$ and $c\in\bbR$:
\[
\bbE\left[e^{-\frac 12 \, (\ga+c)^2}\right] \, = \, \frac
1{\sqrt{1+\sigma^2}} \, e^{-\frac 12 \, \frac{c^2}{1+\sigma^2}},
\]
we have now for $\gga(X)\sim N(0,I)$:
\begin{align*}
& \bbE\left[ e^{-\frac 12 \, \frac 1{1-I} \, \left(\gga(X)+\langle
\overline h,\Lambda\rangle-\kappa\right)^2} \right] \, = \,
\frac1{\sqrt{1+\frac{I}{1-I}}} \, e^{-\frac 12 \, \frac 1{1-I} \,
\frac1{1+\frac{I}{1-I}}\left(\langle \overline
h,\Lambda\rangle-\kappa\right)^2} \, = \, {\sqrt{1-I}} \, e^{-\frac
12 \, \left(\langle \overline h,\Lambda\rangle-\kappa\right)^2}.
\end{align*}
Therefore, recalling the definition of $\overline
h:=h-\frac1{1-I}\langle M,h\rangle(\gl+\mu)$, we obtain after some
trivial computation:
\begin{align*} &
\log \bbE\left[e^{\langle Z,h\rangle} \, \rho(Z) \right] \, = \,
\frac\kappa{1-I}\langle M,h\rangle+\frac 12\langle Q\overline
h,\overline h\rangle -\frac 12 \, \left(\langle \overline
h,\Lambda\rangle-\kappa\right)^2 + \frac 12 \, \kappa^2
\\ & = \kappa\langle \Lambda+M,h\rangle + \frac 12 \langle Qh,h\rangle
- \langle \Lambda+M,h\rangle^2 \, = \, \kappa\langle
h,\Lambda+M\rangle+\frac 12\langle Q_Yh,h\rangle. \qed
\end{align*}


\begin{thebibliography}{abc99xys}


\bibitem{ags}{ L. Ambrosio, N. Gigli, G. Savar\'e (2005),}
\emph{Gradient flows in metric spaces and in the spaces of
probability measures}. Lectures in Mathematics ETH Z\"urich,
Birkh\"auser Verlag, Basel.

\bibitem{asz} L. Ambrosio, G. Savar\'e, L. Zambotti (2007),
{\it Existence and Stability for Fokker-Planck equations with
log-concave reference measure}, preprint posted on
http://arxiv.org/abs/0704.2458.

\bibitem{cepa} E. C\'epa, (1998), {\it Probl\`eme de Skorohod multivoque},
Ann. Prob. {\bf 28}, no. 2, 500-532.

\bibitem{dpr} G. Da Prato, M. R\"ockner (2002),
{\it Singular dissipative stochastic equations in Hilbert spaces},
 Probab. Theory Related Fields  {\bf 124},  no. 2, 261--303.

\bibitem{dpz3} G. Da Prato, J. Zabczyk (2002),
{\it Second order partial differential equations in Hilbert spaces},
London Mathematical Society Lecture Note Series, n. 293.

\bibitem{deza} A. Debussche, L.Zambotti (2007),
{\it Conservative Stochastic Cahn-Hilliard equation with
reflection}, to appear in Annals of Probability.

\bibitem{dim} R.T. Durrett, D.L. Iglehart, D.R. Miller (1977),
{\it Weak convergence to Brownian meander and Brownian excursion},
Ann. Probability, {\bf 5}, no. 1, pp. 117-129.

\bibitem{nupa} D. Nualart, E. Pardoux (1992), {\it
White noise driven quasilinear SPDEs with reflection},
Prob. Theory and Rel. Fields, {\bf 93}, pp. 77-89.

\bibitem{reyo}  D. Revuz, M. Yor,  (1991),
\textit{Continuous Martingales and Brownian Motion}, Springer Verlag.

\bibitem{za} L. Zambotti (2002), {\it Integration by parts
on convex sets of paths and applications to SPDEs with
reflection}, Prob. Theory and Rel. Fields {\bf 123},
no. 4, 579-600.

\bibitem{za4} L. Zambotti, (2004), {\it Occupation densities for SPDEs
with reflection}, Annals of Probability, {\bf 32} no. 1A, 191-215.

\bibitem{f} L. Zambotti (2007), {\it Fluctuations for a
conservative interface model on a wall}, available on
http://arxiv.org/abs/0711.0583.


\end{thebibliography}
\end{document}